\documentclass{amsart}

\usepackage{amsthm,amsfonts,amsmath,amssymb,latexsym,epsfig}
\usepackage{upref,amssymb,eucal,ae,enumitem}
\usepackage[all,cmtip]{xy}

\newtheorem{theorem}{Theorem}
\newtheorem{lemma}[theorem]{Lemma}

\newtheorem{corollary}[theorem]{Corollary}

\newcounter{spthe}

\def\<{\langle}
\def\>{\rangle}

\def\tn{\tilde{n}}
\def\tg{\tilde{g}}
\def\mb#1{{\mathbb #1}}

\begin{document}

\title[On some manifolds with positive sigma invariants]
{On some manifolds with positive sigma invariants and their realizing 
conformal classes}
\author{Santiago R. Simanca}
\email{srsimanca@gmail.com}

\begin{abstract}
We prove that the metric of the Riemannian product 
$(\mb{S}^k(r_1)\times \mb{S}^{n-k}(r_2), g^n_k)$, $r_1^2+r_2^2=1$, is a 
Yamabe metric in its conformal class if, and only if, either $g^n_k$ is 
Einstein, or the linear isometric embedding of this manifold into the 
standard $n+1$ dimensional sphere is minimal. We combine this result 
with Simons' gap theorem to show that, 
for $2\leq k\leq n-2$, the conformal class of the product metric with
minimal embedding, which is at the upper end of Simons' gap theorem, 
realizes the sigma invariant of $\mb{S}^k\times \mb{S}^{n-k}$, 
and that this is the only class that achieves such a value. Similarly,
we use coherent minimal isometric embeddings of suitably scaled standard 
Einstein metrics $g$ on $\mb{P}^n(\mb{R})$, $\mb{P}^n(\mb{C})$, and 
$\mb{P}^n(\mb{H})$ into unit spheres, and determine the sigma 
invariant of these projective spaces, prove that in each case the conformal 
class $[g]$ realizes it, and that this realizing class is unique.
\end{abstract}


\maketitle

\section{A brief preliminary}
By the Nash isometric embedding theorem \cite{nash}, any 
Riemannian $n$-manifold $(M^n,g)$ can be isometrically embedded into 
a standard sphere $(\mb{S}^{\tn}, \tg) \hookrightarrow  
(\mb{R}^{\tn+1}, \| \phantom{.} \|^2)$ in Euclidean space
of sufficiently large but fixed dimension $\tn= \tn(n)$. We shall use this 
fact throughout our work. We shall denote the volume of $(\mb{S}^n,\tg)$ 
by $\omega_{n}$. 

If $f_g: (M^n,g) \rightarrow (\mb{S}^{\tn}, \tg)$ is an isometric embedding  
of a Riemannian manifold $(M,g)$ into $(\mb{S}^{\tn}, \tg)$, then we have the
relation 
\begin{equation} \label{gra}
s_g  = n(n-1) +\tg(H_{f_g},H_{f_g})- \tg(\alpha_{f_g},\alpha_{f_g}) \, , 
\end{equation}
where $s_g$ is the scalar curvature of $g$, and $n(n-1)$,  
$H_{f_g}$, and $\alpha_{f_g}$ are the exterior scalar curvature, mean curvature
vector, and second fundamental form of $f_g(M)$, respectively
\cite[(4)]{gracie}.

If we have given a path $f_{g_t}$ of conformal deformations of $f_{g}$, 
by the Palais isotopic extension theorem, there exists a path $F_t$ of 
diffeomorphism of $\mb{S}^{\tn}$ such that $F_t(f_g(x))=f_{g_t}(x)$, and since
the metrics on the submanifolds are all induced by the metric $\tg$ 
on the background sphere, by pull-back of the metrics by $F_t$, and 
restriction of the diffeomorphisms to $f_g(M)$, we are able to relate $f_{g_t}$
to $f_g$, and express the intrinsic and extrinsic quantities of $f_{g_t}$ in 
terms of the said quantities for $f_g$ itself, and suitable differential 
operators acting on a function $u(t)$, defined in a tubular neighborhood of 
$f_{g}(M) \hookrightarrow \mb{S}^{\tn}$, such that 
$\tilde{g}\mid_{f_{g_t}(M)}=e^{2u(t)}\tilde{g}\mid_{f_g(M)}$.
In particular, we obtain that  
\begin{equation}
n(n-1) = e^{-2u(t)}(n(n-1)-2(n-1)(
{\rm div}_{f_g(M)} (\nabla^{\tg} u)^\tau- \tg(H_{f_g}, \nabla^{\tg}
u^\nu)-\|\nabla^{\tg}u^\tau\|_{\tg}^2+\frac{n}{2}\|\nabla^{\tg} 
u\|_{\tg}^2 )) \label{eq2} \, , 
\end{equation}
\begin{eqnarray}
\| H_{f_{g_t}}\|^2 & = & e^{-2u(t)}(\| H_{f_g}\|^2_{f_g} 
- 2n\tg(H_{f_g},\nabla^{\tg}u ^\nu)
+n^2 \tg(\nabla^{\tg}u ^\nu, \nabla^{\tg} u^\nu)) \, , \label{eq3}\\
\| \alpha_{f_{g_t}}\|^2 & = & e^{-2u(t)}(\| \alpha_{f_g}\|^2 - 2\tg(H_{f_g},
\nabla^{\tg}u^\nu) +n \tg(\nabla^{\tg}u ^\nu, \nabla^{\tg} u^\nu)) \, , 
\label{eq4} 
\end{eqnarray}
where $\nabla^{\tg}u^{\tau}$ and $\nabla^{\tg}u^{\nu}$ stand for the 
tangential and normal component of the gradient
$\nabla^{\tg}u$ on points of $f_g(M)$, expressions that are fully 
determined once we know the first jet of $u(t)$ in the normal directions of  
$f_g(M)$ inside $\mb{S}^{\tn}$. 
By
(\ref{gra}), these expressions imply the intrinsic scalar curvature relation
\begin{equation} \label{eq5}
s_{g_t}= e^{-2u(t)}\left( s_{g}-2(n-1){\rm div}_{f_g(M)}(\nabla^{\tg}u)^\tau
-(n-1)(n-2)\tg(\nabla^{\tg}u^{\tau},\nabla^{\tg}u^\tau)  \right) \, .
\end{equation}
We refer the reader to \S 3 of \cite{scs} for
details, or to \cite[p. 8]{scs2} for a self-contained summary of them.

We set $N=2n/(n-2)$. We recall that the quantity   
$$ 
\lambda(M,[g])= \inf_{g\in [g]} \lambda(M,g)=
\inf_{g\in [g]} \frac{1}{\mu_g(M)^{\frac{2}{N}}}\int s_g d\mu_g
$$ 
is a conformal invariant, that a Yamabe metric in $[g]$ is a metric 
in the class that realizes this invariant, and that any conformal class of 
metrics on $M$ carries a Yamabe representative \cite{ya,au,tr,sc}. 
The resulting function $g \rightarrow \lambda(M,[g])$ is continuous
\cite[Proposition 7.2]{bb}, hence so is $[g] \rightarrow \lambda(M,[g])$.
By a  crucial result of Aubin \cite{au}, 
\begin{equation} \label{aub}
\lambda(M,[g])\leq \lambda(\mb{S}^n,\tg)=n(n-1) \omega_n^{\frac{2}{n}}\, ,
\end{equation}
and, therefore,   
\begin{equation} \label{sig}
\sigma(M)= \sup_{[g]} \lambda(M,[g]) 
\end{equation}
is a well-defined differentiable invariant of $M$ \cite{sc2}.

\section{Products of spheres} 
If $r_1, r_2 \in \mb{R}^{+}$ are such that $r_1^2+r_2^2=1$, and  
$n, k\in \mb{N}$, $n>k$, we consider the product Riemannian manifold
$(M^{n}_k(r_1,r_2),g^n_k)$ where $M^{n}_k(r_1,r_2)=\mb{S}^k(r_1) \times 
\mb{S}^{n-k}(r_2)$, and $g^n_k$ is the product metric of the
sphere factors. We look at its linear isometric embedding
\begin{equation} \label{lem}
\begin{array}{rcl}
f_{g^n_k}: (\mb{S}^k(r_1) \times \mb{S}^{n-k}(r_2),g^n_k) & \rightarrow &  
(\mb{S}^{n+1}, \tg) \hookrightarrow 
(\mb{S}^{\tn(n)}, \tg) \\  (x',x'') \phantom{MMM} & \mapsto & \phantom{i} 
(x',x'') 
\end{array}
\end{equation}
as the hypersurface $f_{g^n_k}(\mb{S}^k(r_1) \times \mb{S}^{n-k}(r_2))$
in the standard sphere of dimension $n+1$. Notice that if $n=2(m+l+1)$ and 
$k=2l+1$, $(M^{n=2(m+l+1)}_{k=2l+1}(r_1,r_2),g^n_k)$ is a particular model of a 
Calabi-Eckmann Hermitian manifold \cite{caec}, and that  
$(M^{n}_k(r_1,r_2),g^n_k)\cong (M^{n}_{n-k}(r_2,r_2),g^n_{n-k})$.
We single out notation wise the case 
$$
\mb{S}^{n,k} = M^n_k\left(\sqrt{\frac{k}{n}},\sqrt{\frac{n-k}{n}}\right)
\subset \mb{S}^{n+1} \, , 
$$
and denote by $g_{\mb{S}^{n,k}}$ its given product metric.

The principal curvatures of $(M^n_k(r_1,r_2),g^n_k)$ are given by  
$\pm r_2/r_1$ and $\mp r_1/r_2$ with multiplicities $k$ and $n-k$, 
respectively. Hence,
$$
\begin{array}{rcl}
\| H_{f_{g^n_k}}\|^2 & = & \left( \pm k\left( \frac{r_2}{r_1}\right) \mp (n-k) 
\left( \frac{r_1}{r_2}\right) \right)^2 =  k^2 \left( \frac{r_2}{r_1}\right)^2 
-2k(n-k) + (n-k)^2 \left( \frac{r_1}{r_2}\right)^2 \, , \\
\| \alpha_{f_{g^n_k}} \|^2 & = & k \left( \frac{r_2}{r_1}\right)^2 
+ (n-k) \left( \frac{r_1}{r_2}\right)^2 \, , \\
s_{g^n_k} & =  & \frac{k(k-1)}{r_1^2}+\frac{(n-k)(n-k-1)}{r_2^2} \, . 
\end{array}
$$
The Ricci tensor of $g^n_k$ is nonnegative, and strictly positive if 
$2\leq k\leq n-2$. 

The manifolds $(\mb{S}^{n,k},g_{\mb{S}^{n,k}})$, $1\leq k\leq n-1$, occur 
at the upper end
of Simons' gap theorem, and their metrics are Yamabe metrics in their  
conformal classes \cite[Theorem 4]{scs} of scalar curvature $n(n-2)$,
with Yamabe invariant 
$$
\lambda(\mb{S}^{n,k},g_{\mb{S}^{n,k}})=
\lambda(\mb{S}^{n,k},[g_{\mb{S}^{n,k}}])=n(n-2)\left( \left(\frac{k}{n}
\right)^{\frac{k}{2}}\left(\frac{n-k}{n}\right)^{\frac{
n-k}{2}}\omega_k \omega_{n-k}\right)^{\frac{2}{n}} \, .
$$
The metric $g^n_k$ on $M^n_k\left(r,\sqrt{1-r^2}\right)$ 
is Einstein if, and only if, $r=\sqrt{\frac{k-1}{n-2}}$, in which case,
it is a Yamabe metric in its conformal class \cite{oba,au} of scalar
curvature $n(n-2)$ also, and 
$$ 
\mu_{g^n_k}=\left( \frac{k-1}{n-2}\right)^{\frac{k}{2}}
\left( \frac{n-k-1}{n-2}\right)^{\frac{n-k}{2}} \omega_k \omega_{n-k}
\leq \mu_{g_{\mb{S}^{n,k}}}\, . 
$$ 
If $2k \neq n$, the function 
$$
[0,1] \ni x \rightarrow \left( \frac{k-x}{n-2x}\right)^{\frac{k}{2}}
\left( \frac{n-k-x}{n-2x}\right)^{\frac{n-k}{2}} 
$$
is stationary at $x=0$, and strictly decreasing on $(0,1]$, so the   
equality of the volumes above occurs if, and only if,
$\sqrt{\frac{k-1}{n-2}}=\sqrt{\frac{k}{n}}$, in which case $n=2k$, and
$g^n_k=g_{\mb{S}^{n,k}}$.
\smallskip

We sharpen these results as follows.

\begin{lemma} \label{le1}
The product metric $g^n_k$ on $M_k^n(r_1,r_2)$ is a Yamabe metric if, and only 
if, either $g^n_k$ is Einstein, or $H_{f_{g^n_k}}=0$ and so 
$(M_k^n(r_1,r_2),g^n_k)=(\mb{S}^{n,k},g_{\mb{S}^{n,k}})$. 
If $2\leq k \leq n-2$, we have that 
$$
\lambda(M^n_k(r_1,r_2),[g^n_k]) \leq 
\lambda(\mb{S}^{n,k},g_{\mb{S}^{n,k}}) \, ,
$$
with equality if, and only if, $[g^n_k]=[g_{\mb{S}^{n,k}}]$.
\end{lemma}

{\it Proof}. By interchanging the factors, if necessary, we may 
assume that $k\leq n/2$.

The variation of the functional $\lambda(M^n_k,g^n_k)$ 
at $g^n_k$ is given by  
\begin{equation} \label{vfYp}
\frac{d}{dt} \lambda(M^n_k,g(t))\mid_{t=0} = 
\frac{1}{\mu_{g^n_k}^{\frac{2}{N}}} \int -(r_{g^n_k},h) d\mu_{g^n_k}
+\frac{1}{\mu_{g^n_k}^{\frac{2}{N}}}
s_{g^n_k}\left(1-\frac{2}{N}\right) \frac{d}{dt} \int d\mu_{g(t)}\mid_{t=0} 
\, , 
\end{equation}
where $g(t)$ is a path of deformations of $g^n_k$ that infinitesimally
varies in the direction of the symmetric two tensor $h=\dot{g}(0)$, and 
$r_{g^n_k}$ is the Ricci tensor of $g^n_k$. This is the 
variation of the Yamabe functional at $g^n_k$ if the path $g(t)$ is conformal.

If $g^n_k$ is Einstein, (\ref{vfYp}) is identically zero for 
any type of metric deformation $g(t)$. If otherwise, we observe that     
the mean curvature vector $H_{f_{g^n_k}}$ is nontrivial if, and only if,
$(M_k^n(r_1,r_2),g^n_k) \neq (\mb{S}^{n,k},g_{\mb{S}^{n,k}})$, and in 
such a case, we prove that there is a choice of conformal path 
$t\rightarrow g(t) \in [g^n_k]$ that makes (\ref{vfYp}) strictly negative, 
which proves then that $\lambda(M^n_k(r,\sqrt{1-r^2}),[g^n_k]) <
\lambda(M^n_k(r,\sqrt{1-r^2}),g^n_k)$, and so $g^n_k$ is not a Yamabe metric
in its class.  

We begin by recalling that if
$$
f_{g(t)}: (M^n_k(r,\sqrt{1-r^2}),g(t)) \rightarrow (\mb{S}^{\tn},\tg)
$$ 
is any path of Nash isometric embeddings of $g(t)$ that deforms 
$f_{g^n_k}$, then 
$$
\frac{d}{dt}d\mu_{g(t)}=\left( {\rm div}(T^{\tau}) - \< T^{\nu},
H_{f_{g_t}}\>\right) d\mu_{g(t)} \, , 
$$
where $T=T^{\tau}+ T^{\nu}$ is the decomposition of the
variational vector field of the path into tangential and normal components.

The nontrivial mean curvature vector 
$H_{f_{g^n_k}}$ may be expressed as 
$H_{f_{g^n_k}}=h_{f_{g^n_k}}\nu_{H_{f_{g^n_k}}}$ for some scalar 
$h_{f_{g^n_k}}>0$ and normal vector $\nu_{H_{f_{g^n_k}}}$  (see the details
of this, and relevant constructions in \cite[p. 15-16]{scs}). 
 At any point on the submanifold, 
we let $t$ be the arc length parameter of the geodesic in 
$(\mb{S}^{\tn},\tg)$ that emanates from  
the said point in the direction of $\nu_{H_{f_{g^n_k}}}$. This defines  
$t$ as a scalar function on the submanifold that ranges in $(-\varepsilon, 
\varepsilon)$ for some sufficiently small $\varepsilon>0$. We then choose
a path of conformally related metrics $(-\varepsilon, \varepsilon)\ni t 
\rightarrow g(t)=e^{2u(t)}\tg\mid_{f_{g^n_k}(M^n_k)}$ 
such that the variational vector field of the Nash isometric embedding 
$f_{g(t)}$ into $(\mb{S}^{\tn},\tg)$ at $t=0$ is $H_{f^n_k}$,
and we have $h=\dot{g}(0)=h_{f_{g^n_k}} \tg\mid_{f_{g^n_k}(M^n_k)}$. 

By the identity above, and the divergence theorem, we have that
$$
\frac{d}{dt}\int d\mu_{g(t)}\mid_{t=0}= 
-\int \| H_{f_{g^n_k}}\|^2 d\mu_{g^n_k} < 0 \, ,  
$$
and since 
$$
s_{g^n_k}\left( 1 -\frac{2}{N}\right) =
s_{g^n_k}\frac{2}{n} > 0 \, , 
$$
along the chosen path, the second summand on the right of (\ref{vfYp}) is 
negative. On the other 
hand, since the background metric $\tg$ induces the intrinsic metric 
$g^n_k$ on the submanifold, we have that
$$
\int (r_{g^n_k},h) d\mu_{g^n_k} = 
\int_{f_{g^n_k}(M^n_k)} (r_{g^n_k},
h_{f_{g^n_k}} g^n_k) 
 d\mu_{g^n_k} = 
s_{g^n_k} \int h_{f_{g^n_k}} d\mu_{g^n_k} > 0 \, ,  
$$
and along this path, the first summand on the right of (\ref{vfYp}) is 
negative as well.

In order to finish the proof, we observe that as $2\leq k$,
the set of $r$s where
$$
\lambda(M_k^n(r,\sqrt{1-r^2}),g^n_k) \leq 
\lambda(\mb{S}^{n,k},[g_{\mb{S}^{n,k}}]) 
$$
is either the point $\sqrt{\frac{k}{n}}=\frac{1}{\sqrt{2}}$ if $2k=n$, or an 
interval of the form $[a, \sqrt{k/n}] \subset (0,1)$ if $2k<n$, and 
correspondingly, the subset where  
$$
\lambda(M_k^n(r,\sqrt{1-r^2}),g^n_k) > 
\lambda(\mb{S}^{n,k},[g_{\mb{S}^{n,k}}]) 
$$
is either $(0,1)\setminus \{ 1/\sqrt{2}\}$, or $(0,a)\cup (\sqrt{k/n},1)$, 
respectively.

We compare with each other the scalar curvatures of the various 
metrics after normalizing their volumes to make them all equal to 
$\mu_{g_{\mb{S}^{n,k}}}$. If we choose a Yamabe metric $g^Y$ of volume 
$\mu_{g_{\mb{S}^{n,k}}}$ in the conformal class of $g^n_k$, and if  
$\Omega=(\omega_k \omega_{n-k}/\mu_{g_{\mb{S}^{n,k}}})^{\frac{2}{n}}$, we 
prove that it is not possible to have
\begin{equation} \label{est}
n(n-2) 
 < s_{g^Y}  < \min\left\{\left(k(k-1) \left(
\frac{1-r^2}{r^2}\right)^{\frac{n-k}{n}} +
 (n-k)(n-k-1) \left(
\frac{r^2}{1-r^2}\right)^{\frac{k}{n}}\right)\Omega , n(n-1) \left( \frac{
\omega_n}{\mu_{g_{\mb{S}^{n,k}}}}\right)^{\frac{2}{n}}\right\}
\, ,
\end{equation}
for any $r \in (0,1)$, which precludes the possibility that
$$
\lambda(\mb{S}^{n,k},[g_{\mb{S}^{n,k}}]) < 
\lambda(M_k^n(r,\sqrt{1-r^2}),g^Y)=\lambda(M_k^n(r,\sqrt{1-r^2}),[g^n_k])\, .  
$$
Notice that we are encoding Aubin's universal bound (\ref{aub}) in the
stated bound on the right side of (\ref{est}).

For notational convenience, we set the scalar curvature of the normalized 
metric $\left(\mu_{g_{\mb{S}^{n,k}}}/\mu_{g^n_k}\right)^{\frac{2}{n}}g_k^n$ on 
the right side of (\ref{est}) to be 
$$
f_k(r)= \left(k(k-1) \left(
\frac{1-r^2}{r^2}\right)^{\frac{n-k}{n}} +
 (n-k)(n-k-1) \left(
\frac{r^2}{1-r^2}\right)^{\frac{k}{n}}\right)\Omega \, .  
$$
Since
$$
f'_k(r)= \frac{2k(n-k)}{r(1-r^2)n} \left((k-1) \left(
\frac{1-r^2}{r^2}\right)^{\frac{n-k}{n}} -
 (n-k-1) \left(
\frac{r^2}{1-r^2}\right)^{\frac{k}{n}}\right)\Omega \, ,   
$$
$f_k(r)$ has only one critical point in $(0,1)$, 
$$
\overline{r}= \sqrt{\frac{k-1}{n-2}} \, , 
$$
where $f_k(r)$ achieves its minimum. Notice that the metric corresponding to
$\overline{r}$ is a volume nondecreasing homothetic transformation of the 
Einstein product metric on $M(\overline{r},\sqrt{1-\overline{r}^2})$, and so
it has scalar curvature less or equal than $n(n-2)$.   
If $\overline{r}< \sqrt{k/n}$, we have that
$$
f_k(\overline{r}) <  f_k(a)=f_k(\sqrt{k/n})=
n(n-2) \, , 
$$
so if there exists an $r \in (0,\sqrt{k/n}) \supset (a,\sqrt{k/n})$ where the 
volume $\mu_{\mb{S}^{n,k}}$ Yamabe metric $g^Y$ satisfies (\ref{est}), by 
continuity, such an $r$ must be strictly smaller than $a$, which  
contradicts the first part of the proof since $\lambda(M^n_k(r,\sqrt{1-r^2},[
g^n_k])$ would have to be strictly less than $f(a)\mu_{g_{\mb{S}^{n,k}}}^{
\frac{2}{n}}$ in a nontrivial neighborhood of 
$a$, so no $\lambda(M^n_k(r,\sqrt{1-r^2},[g^n_k])$ can rise above this
value for $r$s in this neighborhood, and so not at all on $(0,a]$, and thus,
on $(0, \sqrt{k/n}]$. Proceeding similarly with the isometric manifold 
$(M^n_{n-k}(\sqrt{1-r^2},r),g^n_{n-k})$ whose corresponding function
$f_{n-k}(r)$ has a minimum critical point at $\sqrt{(n-k-1)/(n-2)}$, we 
conclude that no $\lambda(M^n_k(r,\sqrt{1-r^2},[g^n_k])$ can rise above 
the value $n(n-2)$ for $r$s in $[\sqrt{(n-k-1)/(n-2)},1)$ either.  
Thus, if (\ref{est}) were to happen in this case, it would have to happen for 
an $r\in (\sqrt{k/n}, \sqrt{(n-k-1)/(n-2)})$, possibility that we need to 
exclude.   
On the other hand, when $\overline{r}=\sqrt{k/n}$, then $2k=n$, and we 
must exclude the possibility that (\ref{est}) holds for any 
$r\in (0,1)\setminus \{\frac{1}{\sqrt{2}}\}$. We proceed to deal with 
these two remaining situations by an general argument that applies to both. 

In $[\sqrt{k/n},1)$, $f_k(r)$ is increasing. 
Suppose that there exists an $r > \sqrt{k/n}$ where (\ref{est}) holds.
By continuity of $f_k(r)$, there must be an $r_0 \in [\sqrt{k/n},r]$
such that $f_k(r_0)=s_{g^Y}$, and so 
$\lambda(M^n_k(r_0,\sqrt{1-r_0^2},g^n_k)=s_{g^Y}\mu_{g_{\mb{S}^{n,
k}}}^{\frac{2}{n}}$. This contradicts the first part of the theorem
since the metric $g^n_k$ corresponding to this $r_0$ is not Einstein, and 
the mean curvature vector of the isometric embedding $f_{g^n_k}(r)$ is 
nonzero at $r=r_0$, so no $\lambda(M^n_k(r,\sqrt{1-r^2},[g^n_k])$ can rise 
to or above this value in a nontrivial neighborhood of $r_0$.
\qed
\medskip

When $k=1$, there is a countable set of $r$s $\nearrow 1$ such that $M^n_1(r,
\sqrt{1-r^2})$ carries a one parameter family of Yamabe metrics of the form 
$e^{2u}g^n_1$ with $e^{2u}$ a nonconstant conformal factor, and  
$\lambda(\mb{S}^{n,1}, g_{\mb{S}^{n,1}})
<\lambda(M^n_1,e^{2u}g^n_1) \nearrow n(n-1) \omega_n^{\frac{2}{n}}$,  
which by (\ref{aub}), allows for the conclusion that
$\sigma(\mb{S}^1 \times 
\mb{S}^{n-1})=n(n-1) \omega_n^{\frac{2}{n}}=
\sigma(\mb{S}^{n-1} \times \mb{S}^1)$ \cite[\S 2]{sc2}. 
In contrast, we have now the following:  

\begin{theorem} \label{th2}
Suppose that $2\leq k \leq n-2$. If $g$ is any Riemannian metric on the 
manifold $\mb{S}^k \times \mb{S}^{n-k}$, then 
$$
\lambda(\mb{S}^k\times \mb{S}^{n-k},[g]) \leq 
\lambda(\mb{S}^{n,k},[g_{\mb{S}^{n,k}}]) =: 
\sigma(\mb{S}^k \times \mb{S}^{n-k}) \, , 
$$
and the equality is achieved if, and only if, $[g]=[g_{\mb{S}^{n,k}}]$.
\end{theorem}

{\it Proof}.
We let $g^Y$ be any smooth Yamabe metric on $\mb{S}^k \times \mb{S}^{n-k}$, and 
scale it, if necessary, so that $\mu_{g^Y}= \mu_{g_{\mb{S}^{n,k}}}$. By  
Lemma \ref{le1} and (\ref{aub}), it suffices to prove that if $[g^Y]$ is not 
a product class
other than possibly $[g_{\mb{S}^{n,k}}]$, and if $s_{g^Y}$ is such that 
\begin{equation} \label{in1}
s_{\mb{S}^{n,k}}= n(n-2) \leq s_{g^Y}  < n(n-1) \left( \frac{\omega_n}{\mu_{
g_{\mb{S}^{n,k}}}} \right)^{\frac{2}{n}} \, , 
\end{equation}
then the Nash isometric embedding of $g^Y$,   
\begin{equation} \label{isp}
f_{g^Y}: (\mb{S}^k \times \mb{S}^{n-k}, g^Y) 
\rightarrow (\mb{S}^{n+p},\tg)
\hookrightarrow (\mb{S}^{\tn},\tg) \, , 
\end{equation}
is, up to isometry of the background, the linear embedding 
$f_{g_{\mb{S}^{n,k}}}$ of the metric  
$g_{\mb{S}^{n,k}}$ in (\ref{lem}), and we in fact have that $s_{g^Y}=n(n-2)$.

We write $H_{f_{g^Y}}= h_{f_{g^Y}} \nu_{f_{g^Y}}$, where 
$h_{f_{g^Y}}$ is a nonnegative constant function, and 
$\nu_{f_{g^Y}}$ is a normal vector field \cite[Theorem 6]{scs}. By the 
orientability of the manifold, these factors of the mean curvature vector 
are globally well-defined, and since $s_{g^Y}$ is constant, by (\ref{gra})
we have that $\| \alpha_{f_{g^Y}}\|^2$ is constant also (cf. 
\cite[Theorem 7]{scs}).

Suppose that $H_{f_{g^Y}} \neq 0$. We let $t$ be once again the arc length 
parameter for the geodesic flow in normal directions. We choose a path 
$t\rightarrow f_{e^{2u(t)}g^Y}$ of homothetics deformations of $f_{g^Y}$ 
defined by a function $u(t)$ such that $u(t)\mid_{f_{g^Y}(M)}=t$, and so the 
tangential gradient of $u(t)$ vanishes, and such that  
$\nabla^{\tg}\dot{u}^\nu \mid_{t=0}=H_{g^Y}$.
Since the Yamabe problem applies for manifolds of dimension $n\geq 3$, the 
resulting path is defined for $t$s in at least $[0,1/n]$, and along such a 
path, by (\ref{eq3}), (\ref{eq4}) and (\ref{eq5}), we have that
$\| H_{f_{g_t}}\|^2$, $\| \alpha_{f_{g_t}}\|^2$, and $s_{g_t}$ remain 
constant functions, and that at $t=1/n$, 
$H_{f_{e^{\frac{2}{n}}g^Y}}=0$, so the isometric embedding  
$f_{e^{\frac{2}{n}}g^Y}$ is minimal. By computing the Yamabe invariant of
$[g^Y]$ using $g^Y$ and $e^{\frac{2}{n}}g^Y$, respectively, we obtain that 
$$
(n(n-1)+\| H_{f_{g^Y}}\|^2 - \| \alpha_{f_{g^Y}}\|^2)
\mu_{g^Y}^{\frac{2}{n}}
= 
(n(n-1)-\| \alpha_{f_{e^{\frac{2}{n}}g^Y}}\|^2)
(e\mu_{g^Y})^{\frac{2}{n}}
\, ,
$$
and by the first of the inequalities in (\ref{in1}), we conclude that
$$
\| \alpha_{f_{e^{\frac{2}{n}}g^Y}}\|^2 \leq n\left( n-1 - e^{-\frac{2}{n}}
(n-2) \right) =(1+a_n)n\, , 
$$
where
$$
a_n := (n-2)\left(1 - e^{-\frac{2}{n}} \right) < 2 \, . 
$$
On the other hand, by the  second inequality in (\ref{in1}), 
we obtain that
$$
e^{\frac{2}{n}}(n(n-1)-\| \alpha_{f_{e^{\frac{2}{n}}g^Y}}\|^2) < 
 n(n-1) \left( \frac{\omega_n}{\mu_{g^Y}} \right)^{\frac{2}{n}} \, ,
$$
and by Lemma \ref{leap1}(a) in Appendix \ref{ap1}, we conclude that 
$$
\| \alpha_{f_{e^{\frac{2}{n}}g^Y}}\|^2 > 
\frac{n(n-1)}{e^{\frac{2}{n}}}\left( e^{\frac{2}{n}} -
 \left( \frac{\omega_n}{\mu_{g^Y}} \right)^{\frac{2}{n}}\right) =(1+b_n)n
\geq (1+a_n)n \, ,
$$
an estimate that contradicts the earlier one for 
the quantity on the left. Thus, $H_{f_{g^Y}}=0$, that is to say, the embedding 
(\ref{isp}) is 
minimal to begin with, and by the first of the inequalities in (\ref{in1}), we 
must have that 
\begin{equation} \label{ces}
\| \alpha_{f_{g^Y}}\|^2 \leq n \, .
\end{equation}

By the Simons' gap theorem \cite[Theorem 5.3.2, Corollary 5.3.2]{si} 
\cite[Main Theorem]{cdck} \cite[Corollary 2]{bla} (cf. with 
\cite[Theorem 9]{scs}), the constant function $\| \alpha_{f_{g^Y}}\|^2$ must
fall into one of two mutually exclusive cases: Either  
$$
\| \alpha_{f_{g^Y}}\|^2 \leq \frac{np}{2p-1} \, ,  
$$
in which case, up to an isometry of the background space, the embedding  
$f_{g^Y}$ coincides with $f_{g_{\mb{S}^{n,k}}}$, $p=1$, and 
$\| \alpha_{f_{g^Y}}\|^2 = n$; or  
$$
\| \alpha_{f_{g^Y}}\|^2 > n > \frac{np}{2p-1} \, .  
$$
By (\ref{ces}), this latter case does not occur.
\qed
\medskip

It has been a question of interest for a long time if a conformal class
on $M$ that attains its sigma invariant carries an Einstein 
representative \cite[pp. 126-127, Lemma 1.2]{sc2}, \cite[\S 4F]{bess},
\cite[Theorem 1.2]{and}, an affirmative version of which became 
known as the ``Besse conjecture.'' Theorem \ref{th2} provides counterexamples
since among the manifolds we considered in there, only the symmetric case of
$\mb{S}^k \times \mb{S}^k$ carries an Einstein metric in the conformal class 
that realizes the said invariant. 

\begin{corollary}
Suppose that $n\neq 2k$, and let $r=\sqrt{\frac{k-1}{n-2}}$. Then the class
$[g_{\mb{S}^{n,k}}]$ does not carry Einstein representatives, and the
Einstein product metric $g^n_k$ on 
$\mb{S}^k(r)\times \mb{S}^{n-k}(\sqrt{1-r^2})$ is such that
$$
\lambda(\mb{S}^k(r)\times \mb{S}^{n-k}(\sqrt{1-r^2}),g^n_k) <
 \sigma( \mb{S}^k \times \mb{S}^{n-k}) \, .  
$$
\end{corollary}

\section{Projective spaces} \label{sec3}
We outline quickly the construction in \cite{si3} extended to include the 
quaternionic projective spaces. We let $\mb{F}$ be any of the division algebra 
fields $\mb{R}$, $\mb{C}$, and $\mb{H}$. Then, $\mb{P}^n(\mb{F})$ is the
quotient space of $\mb{F}^{n+1}\setminus \{0\}$ under the identification
 making $v$ and $v'$ 
equivalent if $v'=\lambda v$ for $\lambda$ in the multiplicative 
group $\mb{F}_1=\{ v \in \mb{F}\setminus \{ 0\}: \; v \overline{v}=
\| v\|^2=1\}$.   
Since $\mb{H}=\mb{C} \oplus \mb{C}j$, we get a natural 
inclusion mapping of $\mb{P}^n(\mb{C})$ in $\mb{P}^n(\mb{H})$, while 
the conjugation operation on $\mb{C}$ renders $\mb{P}^n(\mb{R})$
as the set of real points of $\mb{P}^n(\mb{C})$. 

We define a sequence of positive numbers $r_n$ by 
\begin{equation}\label{ra}
r_n^4=\left(\frac{n+1}{2}\right)^2(n-1)!\, ,
\end{equation}
and view the three alluded projective spaces with the metrics that make
the fibrations of the diagram    
\begin{equation} \label{fib}
\begin{array}{ccccccc}
\mb{Z}/2 & & \phantom{SS} \mb{S}^1 & & \phantom{SS}\mb{S}^3 & & \\ 
& \searrow & & \searrow & & \searrow & \\
& & \mb{S}^n(r_n) & \hookrightarrow & \mb{S}^{2n+1}(r_n) & \hookrightarrow & 
\mb{S}^{4n+3}(r_n) \\
& & \downarrow & & \downarrow & & \downarrow \vspace{1mm} \\
& & \mb{P}^n(\mb{R})
 & \hookrightarrow & \mb{P}^n(\mb{C}) & \hookrightarrow & \mb{P}^n(\mb{H})  
\end{array}
\end{equation}
Riemannian submersions each. For $n=1$, the $\mb{F}_1$ invariant 
mapping
\begin{equation} \label{m1c} 
\mb{S}^{{\rm dim}_{\mb{R}}(\mb{F})+{\rm dim}_{\mb{R}}(\mb{F}_1)} \ni 
v=[v_0:v_1] \mapsto 
\iota^{\mb{F}}_1(v)= (2v_0 \overline{v}_1,(v_0\overline{v}_0-
v_1 \overline{v}_1)) \in \mb{S}^{{\rm dim}_{\mb{R}}\mb{F}} 
\hookrightarrow \mb{R}^{{\rm dim}_{\mb{R}}\mb{F}+1}
\end{equation}
descends to a minimal codimension zero isometric embedding identification
$$
\mb{P}^1(\mb{F}) \hookrightarrow \mb{S}^{{\rm dim}_{\mb{R}}\mb{F}}
 \hookrightarrow \mb{R}^{{\rm dim}_{\mb{R}}\mb{F} +1}
$$
between the domain projective space and the range sphere with its standard
metric, and by construction, we obtain a tower of inclusions 
$$
\begin{array}{ccccc}
\mb{S}^7 & & & & \\
 & \searrow & & & \\
\mb{S}^3 & & \mb{P}^1(\mb{H})  & \hookrightarrow & \mb{S}^4 \\
& \searrow & \cup &  & \cup \phantom{1}\\
\mb{S}^1 & & \mb{P}^1(\mb{C})  & \hookrightarrow & \mb{S}^2 \\
& \searrow & \cup &  & \cup \phantom{1}\\
& & \mb{P}^1(\mb{R})  & \hookrightarrow & \mb{S}^1 
\end{array} \, . 
$$ 
In the intrinsic metrics in these identifications, 
the length, area and 
volume of $\mb{P}^1(\mb{R})$, $\mb{P}^1(\mb{C})$, and $\mb{P}^1(\mb{H})$ are 
$\pi=\frac{1}{2}\omega_1$, $\pi=\frac{1}{2^2}\omega_2$ and $\frac{\pi^2}{3!}=
\frac{1}{2^4} \omega_4$, respectively, while that of the embedded spaces
are, correspondingly, $\omega_1$, $\omega_2$ and $\omega_4$. 

We proceed by induction. We assume that we have defined a map  
$$
\iota^{\mb{F}}_{n-1}: 
\mb{S}^{{\rm dim}_{\mb{R}}\mb{F}^{n-1}+{\rm dim}_{\mb{R}}\mb{F}_1}(r_{n-1}) 
\rightarrow \mb{S}^{L^{\mb{F}}_{n-1}} 
\subset \mb{R}^{L^{\mb{F}}_{n-1}+1}
$$ 
that is invariant under the action of $\mb{F}_1$ on its domain,
descends to a minimal isometric embedding 
$$
\mb{P}^{n-1}(\mb{F}) \hookrightarrow \mb{S}^{L^{\mb{F}}_{n-1}} \hookrightarrow
\mb{R}^{L^{\mb{F}}_{n-1}+1}
$$
of the quotient $\mb{P}^{n-1}(\mb{F})$ with the induced metric, and is such 
that 
$$
\| \iota_{n-1} (v)\|^2 = \frac{1}{r^4_{n-1}}(|v_0|^2+ \cdots + |v_{n-1}|^2)^2 
\, .  
$$
For convenience, we set
\begin{equation} \label{cons}
b^2=\frac{1}{(n^2-1)r_{n-1}^4}\, , \quad a^2=2n(n+1)b^2 \, . 
\end{equation}
If $v=(v',v_n)$, where $v'=(v_0, \ldots, v_{n-1})$, the map  
\begin{equation} \label{mnc}
\iota^{\mb{F}}_n: 
\mb{S}^{{\rm dim}_{\mb{R}}\mb{F}^n+{\rm dim}_{\mb{R}}\mb{F}_1}(r_n) 
\rightarrow 
\mb{S}^{L^{\mb{F}}_n:=L_{n-1}^{\mb{F}}+n{\rm dim}_{\mb{R}}\mb{F}+1} \subset 
\mb{R}^{L^{\mb{F}}_n +1}
\end{equation}
given by 
\begin{equation} \label{mn}
\iota^{\mb{F}}_n(v)=\frac{1}{\sqrt{n+1}}
(\iota^{\mb{F}}_{n-1}(v'),a\overline{v}_nv_0, \cdots , a\overline{v}_nv_{n-1},
b(|v_0|^2+\cdots + |v_{n-1}|^2 - n|v_n|^2))\in \mb{R}^{L_n^{\mb{F}}+1}\, , 
\end{equation}
is such that 
$$
\| \iota^{\mb{F}}_{n} (v)\|^2 = 
\frac{1}{r^4_n}(|v_0|^2+ \cdots + |v_n|^2)^2 \, ,   
$$
and with the Einstein metric $g$ on $\mb{P}^n(\mb{F})$ induced by that of 
$\mb{S}^{{\rm dim}_{\mb{R}}\mb{F}^n+{\rm dim}_{\mb{R}}\mb{F}_1}(r_n)$, 
$\iota^{\mb{F}}_n(v)$ descends to a minimal isometric embedding 
\begin{equation} \label{mef}
\iota^{\mb{F}}_n : \mb{P}^n(\mb{F}) \hookrightarrow \mb{S}^{L^{\mb{F}}_n} 
\subset \mb{R}^{L^{\mb{F}}_n+1}
\end{equation}
into the standard sphere $\mb{S}^{L^{\mb{F}}_n}$. 

We denote the submanifold
$\iota_n^{\mb{F}}(\mb{P}^n(\mb{F}))$ by 
$\mb{P}^n_{\iota_n^{\mb{F}}}(\mb{F})$, and its intrinsic metric by  
$g_{\iota_n^{\mb{F}}}$, respectively.  
The results of \cite{si3} are now enhanced, and subsumed, into the following
theorem, where we in addition specify the geometric intrinsic and extrinsic
quantities of the embedded spaces.

\begin{theorem} \label{th4}
If $r_n$ is the sequence {\rm (\ref{ra})}, and
$$
L_n^{\mb{R}}=\frac{1}{2}n(n+3)-1 \, , \; L_n^{\mb{C}}=(n+1)^2 -2 \, , \; 
L_n^{\mb{H}}=(n+1)(2n+1)-2 \, ,
$$
respectively, then the map {\rm (\ref{mnc})} defined  
inductively by {\rm (\ref{m1c}), (\ref{mn})} above, 
maps the fibers of the fibrations {\rm (}\ref{fib}{\rm )} injectively into 
the image, and with the Einstein metric on $\mb{P}^n(\mb{F})$ 
induced by the metric on the sphere $\mb{S}^{{\rm dim}_{\mb{R}}\mb{F}^n+
{\rm dim}_{\mb{R}}\mb{F}_1}(r_n)$, the map descends to an isometric minimal 
embedding
{\rm (\ref{mef})}, which restricts to the set of real and complex points
making the diagram
\begin{equation} \label{eq21} 
\begin{array}{ccccccccc}
\mb{Z}/2 & & \phantom{SS} \mb{S}^1 & & \phantom{SS}\mb{S}^3 & & & \\ 
& \searrow & & \searrow & & \searrow & & \\
& & \mb{S}^n(r_n) & \hookrightarrow & \mb{S}^{2n+1}(r_n) & \hookrightarrow & 
\mb{S}^{4n+3}(r_n) & \\
& & \downarrow & & \downarrow & & \downarrow & \vspace{1mm} \\
& & \mb{P}^n(\mb{R})
 & \hookrightarrow & \mb{P}^n(\mb{C}) & \hookrightarrow & \mb{P}^n(\mb{H}) \\  
 & & & \stackrel{\iota_n^{\mb{R}}}{\searrow} & & \stackrel{\iota_n^{\mb{C}}}{
\searrow}  & & \stackrel{\iota^{\mb{H}}_n}{\searrow} &
 \vspace{1mm} \\   
 & & & & \mb{S}^{L_n^{\mb{R}}} & \subset & \mb{S}^{L_n^{\mb{C}}} & \subset
 & \mb{S}^{L_n^{\mb{H}}}  
\end{array}
\end{equation}
commutative. If the underlying embedded projective space is of 
real dimension at least two, the intrinsic real sectional 
curvature of its metric is
$$
K^{n,\mb{R}}= 
\frac{1}{2^{\frac{2}{n}}r_n^2}=\frac{1}{2^{\frac{2}{n}}\frac{(n+1)}{2}
\sqrt{(n-1)!}} \, ,
$$
and when it is defined, the intrinsic
holomorphic sectional curvature of the metric in the embedded projective 
space is  
$$
K^{n,\mb{C}}=4 K^{n,\mb{R}}\, .
$$ 
The volume, squared norm of the second fundamental form, 
and scalar curvature of the embedded spaces are 
$$
\begin{array}{rcl}
\mu_{g_{\iota_n^{\mb{R}}}}(\mb{P}^n_{\iota_n^{\mb{R}}}(\mb{R})) & = & 
\omega_n r_n^n=  
2\frac{\pi^{\frac{n+1}{2}}}{\Gamma\left(\frac{n+1}{2}\right)} \left(
\left( \frac{n+1}{2}\right)\sqrt{(n-1)!}\right)^{\frac{n}{2}}\, ,\vspace{1mm}\\
\| \alpha_{\iota_n^{\mb{R}}} \|^2 & = & n(n-1)\left(1-  
1/\left( 2^{\frac{2}{n}}\left(\frac{n+1}{2}\right) 
\sqrt{(n-1)!}\right)\right)\, , \\
s_{g_{\iota_n^{\mb{R}}}} & = &  n(n-1)/\left( 2^{\frac{2}{n}}\frac{n+1}{2}
\sqrt{(n-1)!}\right) \, , 
\end{array}
$$
$$
\begin{array}{rcl}
\mu_{g_{\iota_n^{\mb{C}}}}(\mb{P}^n_{\iota_n^{\mb{C}}}(\mb{C})) & = & 
 (2^{\frac{2}{n}}\left( \frac{n+1}{2}\right)
\sqrt{(n-1)!})^n \frac{\pi^n}{ n!}\, ,  \vspace{1mm} \\ 
\| \alpha_{\iota_n^{\mb{C}}} \|^2 & = & 2n(2n-1)\left(1- \frac{2(n+1)}{2n-1}/
 2^{\frac{2}{n}}\left( \frac{n+1}{2}\right)\sqrt{(n-1)!}\right) \, , \\
s_{g_{\iota_n^{\mb{C}}}} & = &  4n(n+1)/(2^{\frac{2}{n}}\left( \frac{n+1}{2}
\right) \sqrt{(n-1)!})\, , 
\end{array}
$$
$$
\begin{array}{rcl}
\mu_{g_{\iota_n^{\mb{H}}}}(\mb{P}^n_{\iota_n^{\mb{H}}}(\mb{H})) & = & 
 \left( 
2^{\frac{2}{n}} \left( \frac{n+1}{2}\right)
\sqrt{(n-1)!}\right)^{2n}\frac{\pi^{2n}}{
(2n+1)!}\, , \\
\| \alpha_{\iota_n^{\mb{H}}} \|^2 & = & 4n(4n-1)\left(1-\frac{4(n+2)}{4n-1}/(
2^{\frac{2}{n}}\left(\frac{n+1}{2}\right)\sqrt{(n-1)!})\right)\, , \\
s_{g_{\iota_n^{\mb{H}}}} & = & 16n(n+2)/(2^{\frac{2}{n}}\left( \frac{n+1}{2}
\right) \sqrt{(n-1)!}) \, , 
\end{array}
$$
respectively.
\end{theorem}

{\it Proof}. If we dilate the metric of the embedded real projective space 
by a factor of $t$ such that $t^{\frac{n}{2}}=2$, its volume changes to the 
volume of the covering sphere $2\omega_n r_n^n$, and its sectional curvature
gets dilated by the factor $1/t$ yielding the sectional curvature 
of the covering sphere, and that of the projective space itself, which are
the same because the covering map is the local diffeomorphism defined by the 
antipodal map. This proves the stated expression for $K^{n,\mb{R}}$, as well 
as that for  $K^{n,\mb{C}}$, which is four times the former.
    
We have that $\mu_{\iota_n^{\mb{R}}}(\mb{P}_{\iota_n^{\mb{R}}}(\mb{R}))=
\omega_n r_n^n$. Since $s_{g_{\iota_n^{\mb{R}}}}=n(n-1)
K^{n,\mb{R}}$, by (\ref{gra}) we
then get the stated expression for $\| \alpha_{\iota_n^{\mb{R}}} \|^2$. 

Since the embedded real projective space is the set of real points of
the embedded complex projective space, the real and holomorphic sectional
curvature of the latter are $K^{n,\mb{R}}$ and
$K^{n,\mb{C}}$, respectively, and we can thus proceed to compute its intrinsic
scalar curvature by counting the number of real and holomorphic 
sections. We obtain $s_{g_{\iota_n^{\mb{C}}}}=2n(K^{n,\mb{C}}+
(2n-2)K^{n,\mb{R}})=4n(n+1)
K^{n,\mb{R}}$. By (\ref{gra}), the stated result for
$\| \alpha_{\iota_n^{\mb{C}}} \|^2$ follows. The result for the volume
follows by a scaling argument, since the volume of the complex projective
space with the standard Fubini-Study metric is $\pi^n/n!$. 

The volume of $\mb{P}^n({\mb{H}})$ with its standard 
metric is $\pi^{2n}/(2n+1)!$, and using that, the stated result for
$\mu_{g_{\iota_n^{\mb{H}}}}(\mb{P}^n_{\iota_n^{\mb{H}}}(\mb{H}))$
follows by a scaling argument. Since $\mb{P}^n({\mb{C}})$ sits in the 
front and back of $\mb{P}^n({\mb{H}})$, whose metric induces on the
submanifolds the standard Fubini-Study metrics they already have,
we may count the number of real and holomorphic sections of 
$\mb{P}^n{\iota_n^{\mb{H}}}({\mb{H}})$ by adding the previously found 
number of them for the embedded $\mb{P}^n({\mb{C}})$s in it, plus the number of 
them across, all times ${\rm dim}_{\mb{R}}(\mb{H})$. We obtain
that $s_{g_{\iota_n^{\mb{H}}}}=16n(n+2)K^{n,\mb{R}}$. By  
(\ref{gra}), the stated result for the extrinsic quantity 
$\| \alpha_{\iota_n^{\mb{H}}} \|^2$ follows. 
\qed
\medskip

The Riemannian manifolds 
$(\mb{P}^n_{\iota_n^{\mb{R}}}(\mb{R}),g_{\iota_n^{ \mb{R}}})$, 
$(\mb{P}^n_{\iota_n^{\mb{C}}}(\mb{C}),g_{\iota_n^{ \mb{C}}})$,  
$(\mb{P}^n_{\iota_n^{\mb{H}}}(\mb{H}),g_{\iota_n^{ \mb{H}}})$ 
are Einstein, and so their metrics 
are Yamabe metrics in their conformal classes \cite{oba,au}
(cf. \cite[Theorem 4]{scs}). Thus, we obtain 
\begin{equation}
\begin{array}{rcl} \label{eq22}
\lambda(\mb{P}^n_{\iota_n^{\mb{R}}}(\mb{R}),[g_{\iota_n^{\mb{R}}}]) & = & 
\frac{n(n-1)}{2^{\frac{2}{n}}}\omega_n^{\frac{2}{n}} \, , \quad n\geq 3\, ,\\
\lambda(\mb{P}^n_{\iota_n^{\mb{C}}}(\mb{C}),[g_{\iota_n^{\mb{C}}}]) & = & 
4n(n+1) \left(\frac{\pi^n}{n!}\right)^{\frac{1}{n}} \, , \quad n\geq 2 \, ,\\
\lambda(\mb{P}^n_{\iota_n^{\mb{H}}}(\mb{H}),[g_{\iota_n^{\mb{H}}}]) & = & 
16n(n+2) \left(\frac{\pi^{2n}}{(2n+1)!}\right)^{\frac{1}{2n}} \, , \quad
n\geq 1\, .
\end{array}
\end{equation}
The scalar curvature of the metrics in these projective 
spaces is positive, and except when $n=1$, or $n=2$ if $\mb{F}=\mb{R}$,
where the embedding $\iota_n^{\mb{F}}$ into $(\mb{S}^{L_n^{\mb{F}}},\tg)$ is 
of codimension $p=0$, or $p=2$, respectively, we have that 
$\| \alpha_{\iota_n^{\mb{R}}}\|^2 > np/(2p-1)$, 
$\| \alpha_{\iota_n^{\mb{C}}}\|^2 > 2np/(2p-1)$, 
$\| \alpha_{\iota_n^{\mb{H}}}\|^2 > 4np/(2p-1)$, the
lower bounds in these expressions the ones that distinguish the Riemannian 
manifolds at the upper and lower end of Simons' 
gap theorem from the rest of those minimally embedded into 
$(\mb{S}^{L_n^{\mb{F}}},\tg)$ (cf. with \cite[Theorem 2]{rss} where
we reinterpreted this theorem in terms of the critical points 
of the squared global $L^2$-norm of the mean curvature functional, 
under deformations of the immersion, of  
constant density $\| H\|^2$, and $\| \alpha\|^2$ 
sufficiently small relative to it). 

\begin{theorem} \label{th5}
If $g$ is any Riemannian metric on the manifold $\mb{P}^n(\mb{F})$,
$\mb{F}=\mb{R},\mb{C}$ or $\mb{H}$, then 
$$
\lambda(\mb{P}^n(\mb{F}),[g]) \leq 
\lambda(\mb{P}^n_{\iota_n^{\mb{F}}}(\mb{F}),[g_{\iota_n^{\mb{F}}}]) =: 
\sigma(\mb{P}^n(\mb{F}))\, ,
$$
the values on the right those given in {\rm (\ref{eq22})}. The equality is 
achieved by $[g]$ if, and only if, $[g]=[g_{\iota_n^{\mb{F}}}]$.
\end{theorem}

{\it Proof}. We proceed as in the proof of Theorem \ref{th2}, and
choose $g^Y$ to be any smooth Yamabe metric on $\mb{P}^n(\mb{F})$, scaled if
necessary so $\mu_{g^Y}(\mb{P}^n({\mb{F}}))=
\mu_{\iota_n^{\mb{F}}}(\mb{P}^n_{\iota_n^{\mb{F}}}(
\mb{F}))$. 
For convenience we set $n'=n\, {\rm dim}_{\mb{R}}\mb{F}$. 
Then if $g^Y$ is such that 
\begin{equation} \label{rpn}
s_{g_{\iota_{n}^{\mb{F}}}} \leq  
 s_{g^Y}  < n'(n'-1) \left( \frac{\omega_{n'}}{\mu_{
g_{\iota_{n}^{\mb{F}}}}}\right)^{\frac{2}{n'}}  \, , 
\end{equation}
we show that its Nash isometric embedding 
\begin{equation} \label{rpni}
f_{g^Y}: (\mb{P}^{n}(\mb{F}), g^Y) 
\rightarrow (\mb{S}^{n'+p},\tg)
\hookrightarrow (\mb{S}^{\tn},\tg) 
\end{equation}
is minimal, and up to isometries of the background sphere, it coincides with
$f_{g_{\iota_{n}^{\mb{F}}}}$, so
$s_{g^Y}=s_{g_{\iota_{n}^{\mb{F}}}}$. 
 
We write $H_{f_{g^Y}}= h_{f_{g^Y}} \nu_{f_{g^Y}}$, where 
$h_{f_{g^Y}}$ is a nonnegative constant function, and 
$\nu_{f_{g^Y}}$ is a normal vector field \cite[Theorem 6]{scs}. This function
and vector are only locally defined in the nonorientable case, but they always 
multiply to the globally defined mean curvature vector, regardless of the 
orientability of the manifold, and we have that
$\|H_{f_{g^Y}}\|^2$ is a constant function.

If we assume that $H_{f_{g^Y}} \neq 0$, we let $s$ be the arc length parameter
for the geodesic flow in normal directions, and choose a path $s\rightarrow 
f_{e^{2u(s)}g^Y}$ of homothetics deformations of $f_{g^Y}$ that is defined by
a function $u$ equal to $s$ on points of the embedded submanifold, 
where we require that $\nabla u^{\nu}=H_{f_{g^Y}}$ as well. Then, 
by (\ref{eq3}), 
$H_{f_{e^{\frac{2}{n'}}g^Y}} =0$, so the isometric embedding  
$f_{e^{\frac{2}{n'}}g^Y}$ is minimal, and by (\ref{gra}), the Yamabe 
invariant of $[g^Y]$ yields the identity
$$
(n'(n'-1)+\| H_{f_{g^Y}}\|^2 - \| \alpha_{f_{g^Y}}\|^2) \mu_{g^Y}^{\frac{2}{n'}}
= (n'(n'-1)-\| \alpha_{f_{e^{\frac{2}{n'}}g^Y}}\|^2)
(e\mu_{g^Y})^{\frac{2}{n'}} \, .      
$$
By the first of the inequalities in (\ref{rpn}), and the geometric values
in Theorem \ref{th4}, we conclude then that
$$
\| \alpha_{f_{e^{\frac{2}{n'}}g^Y}}\|^2 \leq 
n'(n'-1)-e^{-\frac{2}{n'}}s_{g_{\iota_{n}^{\mb{F}}}}=(1+a_n^{\mb{F}})
\| \alpha_{\iota_{n}^{\mb{F}}}\|^2  \, , 
$$
where
$$
a_n^{\mb{F}}:= 
\frac{(n'(n'-1)- \| \alpha_{\iota_{n}^{\mb{F}}}\|^2)}{
\| \alpha_{\iota_{n}^{\mb{F}}}\|^2}
\left(1 - e^{-\frac{2}{n'}} \right) \, . 
$$
On the other hand, by the second inequality in (\ref{rpn}), we obtain that
$$
e^{\frac{2}{n'}}(n'(n'-1)-\| \alpha_{f_{e^{\frac{2}{n'}}g^Y}}\|^2) < 
 n'(n'-1) \left( \frac{\omega_{n'}}{\mu_{g^Y}} \right)^{\frac{2}{n'}} \, ,
$$
and by Lemma \ref{leap1}(b) in Appendix \ref{ap1}, we conclude that 
$$
\| \alpha_{f_{e^{\frac{2}{n'}}g^Y}}\|^2 > 
\frac{n'(n'-1)}{e^{\frac{2}{n'}}}\left( e^{\frac{2}{n'}} -
\left( \frac{\omega_{n'}}{\mu_{g^Y}} \right)^{\frac{2}{n'}}\right)
=(1+b_n^{\mb{F}}) 
\| \alpha_{\iota_{n}^{\mb{F}}}\|^2 \geq 
(1+a_n^{\mb{F}}) \| \alpha_{\iota_{n}^{\mb{F}}}\|^2  \, , 
$$
which contradicts the previously derived estimate for 
$\| \alpha_{f_{e^{\frac{2}{n'}}g^Y}}\|^2 $. 
Thus, $H_{f_{g^Y}}=0$, that is to say, the embedding (\ref{rpni}) is 
minimal, and by (\ref{gra}), we must have that 
\begin{equation}\label{eq25}
s_{g_n^{\mb{F}}} =n'(n'-1)-\| \alpha_{\iota_{g_n^{\mb{F}}}} \|^2 \leq
n'(n'-1)-\| \alpha_{f_{g^Y}}\|^2 =  s_{f_{g^Y}}\, . 
\end{equation}

In the case $\mb{F}=\mb{R}$, and so $n\geq 3$, we consider the 2-to-1 cover 
space $\mb{S}^n$ of $\mb{P}^n(\mb{R})$ in (\ref{eq21}), and lift the metrics 
$g^Y$ and
$g_{{\iota}_n^{\mb{R}}}$ on the projective space to metrics 
$\tilde{g}^Y$ and $\tilde{g}_n^{\mb{R}}$ on the cover,  
each of volume $2\omega_n r_n^n$. Since the cover map is a 2-to-1 local 
diffeomorphism, the scalar curvatures of the lifted metrics coincide with
the scalar curvatures of the metrics themselves, and the value of the Yamabe
functional on $\tilde{g}^Y$ is greater or equal than Aubin's universal bound
(\ref{aub}). If $[\tilde{g}^Y]\neq [\tilde{g}_n^{\mb{R}}]$, then $\tilde{g}^Y$ 
is not a Yamabe metric in its class, and there exists a volume
preserving conformal deformation changing it to one, 
which must therefore be of scalar curvature strictly smaller than the 
scalar curvature of $\tilde{g}^Y$, and whose projection back to the base of 
the cover shows the existence of a constant scalar curvature representative of 
$[g^Y]$, for which the value of the Yamabe functional is smaller than the 
value of this functional on $g^Y$ itself, contradicting the fact that $g^Y$ 
is a Yamabe metric in $[g^Y]$. Hence, we must have that 
$[\tilde{g}^Y]=[\tilde{g}_n^{\mb{R}}]$, 
and since $g^Y$ is a Yamabe metric in $[g_n^{\mb{R}}]$, 
by the solution of the Yamabe problem on the standard sphere,  
$\tilde{g}^Y$ must be a conformal diffeomorphism deformation of the standard 
metric on the sphere, and the sectional curvature of $\tilde{g}^Y$ must coincide
with that of $\tilde{g}^{\mb{R}}_n$, the constant $K^{n,\mb{R}}$. By the 
local diffeomorphism property of the 
covering map, we then see that the sectional curvature of 
$(\mb{P}^n(\mb{R}), g^Y)$ is the constant $K^{n,\mb{R}}$ also. 
Thus, up to isometries of the background $(\mb{S}^{n(n+3)/2 -1},\tg)
\hookrightarrow (\mb{S}^{\tn},\tg)$,  
$g^Y=g_{\iota_n^{\mb{R}}}$, and
$f_{g^Y}$ is  
the isometric embedding 
$\iota_{n}^{\mb{R}} : (\mb{P}^{n}(\mb{R}), 
g_{\iota_{n}^{\mb{R}}}) \rightarrow 
 \mb{S}^{L_{n}^{\mb{R}}} \hookrightarrow \mb{R}^{L_{n}^{\mb{R}}+1}$ 
of Theorem \ref{th4}, as desired.

In the remaining cases, we assume first that $[g^Y]\neq [g_{\iota_n^{\mb{F}
}}]$. We consider the Riemmanian submersion $(\mb{S}^{n' +{\rm dim}_{\mb{R}}
\mb{F}_1}, \tilde{g}_{\iota^{\mb{F}}_n})
\rightarrow (\mb{P}^n(\mb{F}), g_{\iota^{\mb{F}}_n})$ 
in (\ref{eq21}), and modify the lift $\tilde{g}^{\mb{F}}_n$ of 
$g_{\iota^{\mb{F}}_n}$ in the horizontal
directions while preserving the geometry of the fibers to obtain 
a metric $\tilde{g}^Y$ such that 
$(\mb{S}^{n' + {\rm dim}_{\mb{R}}\mb{F}_1},\tilde{g}^Y) \rightarrow 
(\mb{P}^n(\mb{F}), g^Y)$ is a Riemannian submersion, cf. 
\cite[Theorem 9.59]{bess}. Thus, we have that 
$\mu_{\tilde{g}^{\mb{F}}_n}=
\mu_{\tilde{g}^Y}$, and by (\ref{eq25}), that $s_{\tilde{g}^Y}\geq 
s_{\tilde{g}^{\mb{F}}_n}$, so the value of the Yamabe functional 
on $\tilde{g}^Y$ is greater or equal than Aubin's universal bound
(\ref{aub}) for the total space sphere of the fibration, which is achieved by
$\tilde{g}^{\mb{F}}_n$, and so
$\tilde{g}^Y$ is not a Yamabe metric in its class $[\tilde{g}^Y]$.
It is therefore possible to produce a volume preserving conformal deformation 
of $\tilde{g}^Y$, which preserves the geometry of the fiber, to a metric on 
which the value of the Yamabe functional is strictly less than that of the 
standard $n'+{\rm dim}_{\mb{R}}\mb{F}_1$ sphere. But then the horizontal
component of such a metric has Yamabe functional value strictly less than the
Yamabe functional value of $g^Y$, which contradicts the fact that $g^Y$ is a 
Yamabe metric in its class. 
Hence, we must have that $[g^Y]=[g_{\iota_n^{\mb{F}}}]$, and the
equality in (\ref{eq25}) holds. In this setting, we then consider a path 
$f_{g_t}: (\mb{P}^n(\mb{F}),g_t) 
\rightarrow (\mb{S}^{\tn},\tg)$ of volume preserving conformal isometric
embedding deformations taking $f_{g_{\iota^{\mb{F}}_{n}}}$ at $t=0$ to 
$f_{g^Y}$ at $t=1$, with $\tg\mid_{f_{g_t}(\mb{P}^n(\mb{F}))}=
e^{2u(t)}\tg\mid_{f_{g_{\iota_n^{\mb{F}}}}(\mb{P}^n(\mb{F}))}$ for a scalar
function $u(t)$, defined on a tubular neighborhood of
$f_{g_{\iota^{\mb{F}}_{n}}}(\mb{P}^n(\mb{F}))$,
$u(t)\mid_{t=0}=0$. Since 
the ends of the path are minimal embeddings of equal scalar curvature and
volume, by (\ref{eq3}) we conclude that $\nabla^{\tg}u(t)^{\nu}\mid_{t=1}=0$, 
and as $\| \alpha_{f_{g^Y}}\|^2 = \| \alpha_{\iota_{g_n^{\mb{F}}}} \|^2 >0$, 
then by (\ref{eq4}) that $u(t)\mid_{t=1}=0$. Hence, up to
isometries of the background $(\mb{S}^{L_n^{\mb{F}}},\tg)
\hookrightarrow (\mb{S}^{\tn},\tg)$, $g^Y=g_{\iota_n^{\mb{F}}}$, and
$f_{g^Y}$ is
the isometric embedding
$\iota_{n}^{\mb{F}} : (\mb{P}^{n}(\mb{F}),
g_{\iota_{n}^{\mb{F}}}) \rightarrow
 \mb{S}^{L_{n}^{\mb{F}}} \hookrightarrow \mb{R}^{L_{n}^{\mb{F}}+1}$
of Theorem \ref{th4}.
\qed
\medskip

The sigma invariants of $\mb{P}^3(\mb{R})$ and $\mb{P}^2(\mb{C})$ in the Theorem
above were previously known, \cite{brne} \& \cite{leb}, respectively, 
but in the former of these two cases, our uniqueness of the realizing class 
seems to be a new addendum.

\section{Levy measure Riemannian manifolds and uniform structure}
The infinite dimensional sphere $\mb{S}^{\infty}$ has the natural uniform 
structure associated to its standard metric $\tg$.
 We may include in this noncompact
space all the Nash isometric images of finite dimensional Riemannian manifolds,
and this produces a nontrivial effect on the homotopy theory of the
manifolds when passing from the category of continuous maps to the category 
of uniformly continuous maps, as studied in \cite{grmi} over Levy families.

A family $\{ M_n,\mu_n\}$ of Borel metric measure spaces with 
normalized measures $\mu_n(M_n)=1$ is said to be Levy if, for
any sequence of Borel sets $A_n \subset M_n$ such that $\liminf_{n\rightarrow
\infty}\mu_n(A_n) > 0$, and for every $\varepsilon > 0$, we have that
$\lim_{n\rightarrow \infty}\mu_n(N_{\varepsilon}(A_n)) =1$.
Here, $N_{\varepsilon}(A)$ is the $\varepsilon$-neighborhood of $A$.

If $(M,g)$ is a Riemannian manifold, we take for $\mu$ the normalized
Riemannian volume element. We set $r(M)=\inf_{e\in \mb{S}(M)}r_g(e,e)$, where
$\mb{S}(M)$ is the unit sphere bundle of $(M,g)$. 

The family $\{ \mb{S}^n,\tg\}$ of standard spheres is Levy, as is any 
nontrivial subsequence of it \cite[Principal Example 1.1]{grmi}.  
 
\begin{lemma} 
The family of Riemannian manifolds $\{\mb{S}^{n,k},g_{\mb{S}^{n,k}}\}$,
$2\leq k \leq n-2$, is Levy.
\end{lemma}

{\it Proof}. We may assume that $2\leq k \leq n/2$. Then we have that
$$
r(\mb{S}^{n,k})=\inf_{e\in \mb{S}(M)}r_g(e,e)=\frac{k-1}{k}n\, ,
$$
which approaches $\infty$ as $n\rightarrow \infty$. The result follows by
\cite[Theorem \S1.2]{grmi},
\qed 
\medskip

The situation above is substantially different for any of the
projective spaces $(\mb{P}^n_{\iota_n^{\mb{F}}}(\mb{F}),g_{\iota_n^{\mb{F}}})$ 
associated with the fields $\mb{F}=\mb{R},\mb{C},\mb{H}$, respectively.
The Ricci tensors for all of them approach zero exponentially fast, as follows
by the explicit expressions for their real and complex sectional, and scalar 
curvatures in Theorem \ref{th4}. 

\begin{theorem}
Consider the sequence of projectives spaces 
$(\mb{P}^n_{\iota_n^{\mb{F}}}(\mb{F}), g_{\iota_n^{\mb{F}}}) \hookrightarrow
(\mb{S}^{L_n^{\mb{F}}},\tg )$,
and let $(\mb{P}^{\infty}_{\iota_{\infty}^{\mb{F}}}(\mb{F}),g_{\iota_{\infty}^{
\mb{F}}})$ denote its uniform limit in $(\mb{S}^\infty, \tg)$.
Then $(\mb{P}^{\infty}_{\iota_{\infty}^{\mb{F}}}(\mb{F}),
g_{\iota_{\infty}^{\mb{F}}})$ is an infinite dimensional Ricci flat 
submanifold of $(\mb{S}^\infty, \tg)$.  
$\mb{P}^{\infty}_{\iota_{\infty}^{\mb{R}}}(\mb{R})$ 
and $\mb{P}^{\infty}_{\iota_{\infty}^{\mb{C}}}(\mb{C})$ 
are Eilengberg-MacLane 
spaces $K(\mb{Z}/2,1)$, and $K(\mb{Z},2)$, respectively, while in rational
homotopy, $\mb{P}^{\infty}_{\iota_{\infty}^{\mb{H}}}(\mb{H})$ 
is an Eilengberg-MacLane space $K(\mb{Q},4)$\footnotemark .
\end{theorem}

\footnotetext{We acknowledge an enlightening comment by D.S. asserting this 
property of $\mb{P}^{\infty}_{\iota_{\infty}^{\mb{H}}}(\mb{H})$,
leading to the correction of 
our earlier statement that it was a $K(\mb{Z},4)$: $\pi_i(
\mb{P}^{\infty}_{\iota_{\infty}^{\mb{H}}}(\mb{H})) \cong \pi_{i-1}(\mb{S}^3)$
for all $i$, hence, this space is topologically far more complicated than 
what we erroneously said.}

\appendix 
\section{Scalar curvature comparisons}
\label{ap1}

By Aubin's universal bound (\ref{aub}), the scalar curvatures of Yamabe 
metrics $g$ and $g'$ of equal volume may be meaningfully compared 
with each other, and if the Nash isometric embeddings of $g$ and 
$e^{\frac{2}{n}}g'$ are minimal, we may use (\ref{gra}) to carry out this 
comparison by comparing $\| \alpha_{f_g}\|^2$ and 
$\| \alpha_{f_{e^{\frac{2}{n}}g'}}\|^2$ instead. This imposes strict bounds
if the compared metrics are of positive scalar curvature.   

\begin{lemma} \label{leap1} 
\begin{enumerate}
\item[{\rm (a)}] Consider the Riemannian manifold product  
$(\mb{S}^{n,k},g_{\mb{S}^{n,k}})$ for $n\geq 3$. Then 
$$
\frac{n(n-1)}{e^{\frac{2}{n}}}\left( e^{\frac{2}{n}} -
 \left( \frac{\omega_n}{\mu_{g_{\mb{S}^{n,k}}}} 
\right)^{\frac{2}{n}}\right) - \| \alpha_{f_{g_{\mb{S}^{n,k}}}}\|^2 :=b_n n 
\, ,
$$
and we have that
$$
b_n \geq (n-2)(1-e^{-\frac{2}{n'}}) \, . 
$$
\item[{\rm (b)}] Consider any of the Riemannian manifolds 
$(\mb{P}^n_{\iota_n^{\mb{F}}} (\mb{F}),g_{\iota_n^{\mb{F}}})$,
$\mb{F}=\mb{R},\mb{C}$ or $\mb{H}$, for $n\geq 2$, and let
$n'=n \, {\rm dim}_{\mb{R}}(\mb{F})$ be its real dimension. Then 
$$
\frac{n'(n'-1)}{e^{\frac{2}{n'}}}\left( e^{\frac{2}{n'}} -
 \left( \frac{\omega_{n'}}{\mu_{g_{\iota_n^{\mb{F}}}}} 
\right)^{\frac{2}{n'}}\right) - \| \alpha_{\iota_{g_n^{\mb{F}}}}\|^2 
:=b_n^{\mb{F}} \| \alpha_{\iota_n^{\mb{F}}}\|^2 \, ,
$$
and we have that
$$
b_n^{\mb{F}}\geq \frac{(n'(n'-1)- \| \alpha_{\iota_{n}^{\mb{F}}}\|^2)}{
\| \alpha_{\iota_{n}^{\mb{F}}}\|^2}
\left(1 - e^{-\frac{2}{n}} \right) \, . 
$$
\end{enumerate}
\end{lemma}

{\it Proof}. (a) The product metrics on the product of spheres at the upper 
end of Simons' gap theorem are all of curvature $n(n-2)$, and
$\| \alpha_{f_{g_{\mb{S}^{n,k}}}}\|^2 =n$.
We have that
$$
b_nn=n\left( n-2 -\frac{n(n-1)}{(4\pi k^k(n-k)^{n-k})^{\frac{1}{n}}}\left( 
 \frac{\Gamma(\frac{k+1}{2})\Gamma(\frac{n-k+1}{2})}{e\Gamma(\frac{n+1}{2})}
 \right)^{\frac{2}{n}}\right) :=
(n-2)\left(1 -\frac{c_{n,k}}{e^{\frac{2}{n}}}\right) \, , 
$$
and this can be checked to be positive for any $n,k$ such that 
$n\geq 3$, and $1\leq k\leq n-1$. In this range, we have that
$$
0<c_{n,k}=\frac{1}{n-2}\left( \frac{\omega_n}{\mu_{g_{\mb{S}^{n,k}}}} 
\right)^{\frac{2}{n}} < 1\, ,
$$
and the last assertion follows.

(b) If $\mb{F}=\mb{R}$, by Theorem \ref{th4}, we have that  
$$
b_n^{\mb{R}} \| \alpha_{\iota_n^{\mb{R}}}\|^2 =
\frac{n(n-1)}{e^{\frac{2}{n}}}\left( e^{\frac{2}{n}} -
 \left( \frac{\omega_{n}}{\mu_{g_{\iota_n^{\mb{R}}}}} 
\right)^{\frac{2}{n}}\right) - \| \alpha_{\iota_{g_n^{\mb{R}}}}\|^2 
=\frac{n(n-1)}{2^{\frac{2}{n}}r_n^2}\left( 1- \left( \frac{2}{e}\right)^{
\frac{2}{n}}\right) \, ,
$$
and 
$$
b_n^{\mb{R}}:= \frac{(n(n-1)- \| \alpha_{\iota_{n}^{\mb{R}}}\|^2)}{
\| \alpha_{\iota_{n}^{\mb{R}}}\|^2}\left(1 - \frac{1}{e^{\frac{2}{n}}}
\frac{2^{\frac{2}{n}}}{n(n-1)} \right) \geq  
\frac{(n(n-1)- \| \alpha_{\iota_{n}^{\mb{R}}}\|^2)}{
\| \alpha_{\iota_{n}^{\mb{R}}}\|^2}\left(1 - \frac{1}{e^{\frac{2}{n}}}
\right) \, . 
$$

If $\mb{F}=\mb{C}$, by Theorem \ref{th4}, we have that
$$
b_n^{\mb{C}} \| \alpha_{\iota_n^{\mb{C}}}\|^2 =
\frac{n'(n'-1)}{e^{\frac{2}{n'}}}\left( e^{\frac{2}{n'}} -
 \left( \frac{\omega_{n'}}{\mu_{g_{\iota_n^{\mb{C}}}}} 
\right)^{\frac{2}{n'}}\right) - \| \alpha_{\iota_{g_n^{\mb{C}}}}\|^2 
=\frac{4n(n+1)}{2^{\frac{2}{n}}r_n^2}\left( 1 - \frac{n-\frac{1}{2}}{n+1}
\left( \frac{2\pi \,\omega_{2n}}{e\, \omega_{2n+1}}\right)^{\frac{1}{n}}\right)
\, ,
$$
and 
$$
b_n^{\mb{C}}:= \frac{(2n(2n-1)- \| \alpha_{\iota_{n}^{\mb{C}}}\|^2)}{
\| \alpha_{\iota_{n}^{\mb{C}}}\|^2}\left(1 - \frac{1}{e^{\frac{2}{n'}}}
\frac{1}{4n(n+1)} \left( \frac{2\pi \omega_{2n}}{\omega_{2n+1}}\right)^{
\frac{2}{n'}}\right) \geq \frac{(2n(2n-1)- \| \alpha_{\iota_{n}^{\mb{C}}}\|^2
)}{\| \alpha_{\iota_{n}^{\mb{C}}}\|^2}\left(1 -\frac{1}{e^{\frac{2}{n'}}}
\right) \, . 
$$

If $\mb{F}=\mb{H}$, by Theorem \ref{th4}, we have that 
$$
b_n^{\mb{H}} \| \alpha_{\iota_n^{\mb{H}}}\|^2 =
\frac{4n(4n-1)}{e^{\frac{2}{4n}}}\left( e^{\frac{2}{4n}} -
 \left( \frac{\omega_{n'}}{\mu_{g_{\iota_n^{\mb{H}}}}} 
\right)^{\frac{2}{n'}}\right) - \| \alpha_{\iota_{g_n^{\mb{H}}}}\|^2 
=\frac{16n(n+2)}{2^{\frac{2}{n}}r_n^2}\left( 1 - \frac{4n-1}{4(n+2)}
\left( \frac{2\pi^2\,\omega_{4n}}{e\, \omega_{4n+3}}\right)^{\frac{1}{2n}}
\right) \, , 
$$
and 
$$
b_n^{\mb{H}}:= \frac{(4n(4n-1)- \| \alpha_{\iota_{n}^{\mb{H}}}\|^2)}{
\| \alpha_{\iota_{n}^{\mb{H}}}\|^2} \left(1 - \frac{1}{e^{\frac{2}{n'}}}
\frac{1}{16n(n+2)} \left( \frac{2\pi^2 \omega_{4n}}{\omega_{4n+3}}\right)^{
\frac{2}{n'}}\right) \geq \frac{(4n(4n-1)- \| \alpha_{
\iota_{n}^{\mb{H}}}\|^2)}{\| \alpha_{\iota_{n}^{\mb{H}}}\|^2}\left(1 - 
\frac{1}{e^{\frac{2}{n'}}} \right) \, . 
$$
\qed

\end{document}